%
%
%

\font\magbf=cmbx10 scaled\magstep2


\def\mbox#1{\leavevmode\hbox{#1}}
\def\to{\mbox{$\longrightarrow$}}	
\def\frac#1#2{{#1\over #2}} 
 


\count35=0
\count36=0
\count38=0


\def\section#1{\advance\count35 by 1
\vskip4ex\noindent{\magbf\number\count35 \ \ \ #1}
\count36=0\count38=0 \hfill\vskip2ex}

\def\subsection{\vskip1ex\advance\count36 by 1
{\noindent{\bf\number\count35.\number\count36}\ }} 

\def\today{\ifcase\month\or January\or February\or March\or April\or May\or
June\or July\or August\or September\or October\or November\or December\fi
\space\number\day, \number\year}

\def\beginth#1{\vskip3ex\advance\count36 by 1 {\noindent{\bf
\number\count35 .\number\count36} \ \ {\bf #1.} }\begingroup\sl} 

\def\endth{\endgroup\hfill\vskip2ex}

\def\defn{\vskip3ex\advance\count36 by 1 {\noindent{\bf \number\count35
.\number\count36} \ \ {\bf Definition.} }}

 \centerline{\magbf  Maps of Surface Groups to Finite Groups}

     \centerline{\magbf  with No Simple Loops in the Kernel}
\vskip2ex

      \centerline{                 by }
\vskip2ex
             \centerline{ {\bf Charles Livingston}}

\vskip10ex

\noindent{\bf  Abstract:}   Let $F_g$ denote the closed orientable surface of
genus $g$.  What is the least order finite group,
$G_g$, for which there is a homomorphism $\psi
\!:\pi_1(F_g)
\rightarrow G_g$ so that no nontrivial simple closed curve on $F_g$ represents an
element in Ker($\psi$)?  For the torus it is easily seen that
$G_1 = Z_2
\times Z_2$ suffices.  We prove here that $G_2$ is a group of order 32 and that an
upper bound for the order of $G_g$ is given by $g^{2g +1}$.  The previously known
upper bound was greater than $2^{g{2^{2g}}}$.

\vskip8ex

     For any compact surface $F$ there exists a finite group $G$ and a homomorphism
$\psi\!: \pi_1(F) \rightarrow G$ such that no nontrivial element in the kernel of
$\psi$ can be represented by a simple closed curve. Such a homomorphism is said to
have {\it nongeometric kernel}. Casson, Gabai, and Skora [5] have each constructed
examples of this (see Section 2 for details). The presence of such examples raises a
variety of questions relating to the characterization of the finite groups that can
occur in this way. This paper addresses the problem of determining the relationship
between the genus of $F$ and the order of $G$. In the case that $F$ is a torus a
complete analysis is straightforward. For instance, the natural projection
$\psi \!:
\pi_1(F) \rightarrow H_1(F;Z_2 ) \cong Z_2 \times Z_2$   has nongeometric kernel.

     Our first result concerns the  genus 2 closed orientable surface, $F_2$.
Casson's construction yields a group of order $2^{38}$.  Skora reduced this order
considerably by producing a group of order $2^9$.  In Section 3 a group of order
$2^5 = 32$,
$G_2$,  is constructed for which there is a homomorphism $\psi_2 \!:\pi_1(F_2)
\rightarrow G_2$ having nongeometric kernel. In Section 4 it is proved that no such
example can be constructed using a group of order less than 32.

     The example in Section 3 is generalized to construct examples for arbitrary
genus surfaces in Section 5. The order of the groups constructed is quite small
compared to previously constructed examples. As the examples directly generalize
the minimal genus 2 example, there is the possibility that they are minimal as well.

\vskip1ex

\noindent{\sl Acknowledgements}  Thanks are due to Allan Edmonds for pointing out
the proof of Theorem 4.2.  The work in Section 5 was motivated by discussions with
Dennis Johnson.

\section{Notation and Conventions}

     Throughout this paper all surfaces will be closed and orientable.  References
to basepoints for the fundamental group of a space are omitted. Since the property
of being in the kernel of a homomorphism depends only on the conjugacy class of an
element, such omissions will not affect the arguments.
    
 By a {\it simple loop} on a surface we mean an embedding of the circle $ S^1$ .
    
 We will say that a homomorphism $\psi\!: \pi_1 (F) \rightarrow G$ has {\it geometric
kernel} if some nontrivial element in the kernel can be represented by a simple
loop. Otherwise $\psi$ has {\it nongeometric kernel}. 

\section{Basic Examples}

     In this section a procedure of Casson is used to construct for each surface
$F$ a finite group $G$ and a surjective homomorphism $\psi\!: \pi_1(F) \rightarrow G$ 
such that
$\psi$ has nongeometric kernel. The orders of the groups involved is computed for
contrast with the examples produced in Section 5. 

The statement that $\psi$ has nongeometric kernel can be reinterpreted in terms of
covering spaces as follows. Corresponding to Ker$(\psi)$ there is a connected
regular covering space $p \!: \tilde{F} \rightarrow F$ with 
$p_*(\pi_1(F)) = $ Ker$(\psi)$. An element in $\pi_1 (F)$ is in Ker$(\psi)$ if and
only if when represented by a closed path, the path lifts to a closed path in
$\tilde{F}$. Hence a simple loop on $F$ represents an element in Ker$(\psi)$ if and
only if it can be lifted to a simple loop in $\tilde{F}$.
Conversely, if $p\!: \tilde{F} \rightarrow F$ is a regular covering space with
the property that no nontrivial simple loop on $F$ lifts to $\tilde{F}$ then the
natural projection $\psi\!: \pi_1 (F) \rightarrow \pi_1(F)/p_*(\pi_1(\tilde{F}))$
has nongeometric kernel.
\vskip2ex 

\noindent{\bf Construction}  Given a surface $F$, construct the covering space
$p \!: {\tilde{F}} \rightarrow F$ corresponding to the kernel of the projection
$\pi_1 (F) \rightarrow H_1(F ; Z_2 )$. Simple nonseparating loops on $F$ represent
generators of $H_1(F ; Z_2 )$ and hence do not lift to $\tilde{F}$ . Nontrivial
separating simple loops do lift, but each preimage on $\tilde{F}$ is nonseparating
on $\tilde{F}$.

	Now construct the covering $q \!: \bar{F} \rightarrow \tilde{F}$			
corresponding to the kernel of the projection $\pi_1 (\tilde{F}) \rightarrow
H_1(\tilde{F} ; Z_2 )$.	As no nonseparating simple loop on $\tilde{F}$ lifts to
$\bar{F}$ it is apparent that no nontrivial simple loop on $F$ lifts to
$\bar{F}$.

It remains to show that the covering $p \circ q  \!: \bar{F} \rightarrow 
F$ is regular; that is, that $p_* \circ  q_*(\pi_1(\bar{F}))$ is normal in
$\pi_1(F)$. Observe that $q_*(\pi_1(\bar{F}$)) is a characteristic subgroup of
$\pi_1(\tilde{F})$ and $p_*(\pi_1(\tilde{F}))$      is a characteristic subgroup of
$\pi_1({F})$. Since a characteristic subgroup of a characteristic subgroup is
characteristic, $p_* \circ  q_*(\pi_1(\bar{F})$)    is characteristic in
$\pi_1({F})$, and is hence normal.

\vskip2ex
\noindent{\bf Order of $\pi_1({F}) / < p_* \circ  q_*(\pi_1(\bar{F})) >$}

The order of this finite group is equal to the degree of the covering $p \circ q$.
Suppose that $F$ is of genus $g$. The Euler characterisitic of
$F$ is $2 - 2g$ .  Since  $\tilde{F}$   is a $2^{2g}$ fold cover of $F$, the Euler
characteristic of  $\tilde{F}$ is $2^{2g} (2-2g)$.    The genus of $\tilde{F}$ is
${1 \over 2}(2 - 2^{2g}( 2-2g)) =  { 1 \over 2}((g-1)2^{2g+1} + 2) = g'$.     The
covering $q \!: \bar{F} \rightarrow \tilde{F}$	is of degree $2^{2g'}$.  
The degree of $p \circ q  \!: \bar{F} \rightarrow  F$ is the product of
these two degrees: $2^{2g'} 2g = 2^{(g-1) 2^{2g+1} +2+2g}$.

\vskip2ex
\noindent{\bf Note}  The construction of Gabai differs considerably from the one
above.  He notes that every simple curve is in the complement of some index three
subgroup of $\pi_1(F)$; nonseparating curves are not in the kernel of some map to
$Z_3$ and separating curves are mapped to a 3--cycle in the third symmetric group,
$S_3$, under some homomorphism and hence map to the complement of an index three
subgroup of $S_3$.  Hence the quotient map $\pi_1(F) \to \pi_1(F)/H$ has
nongeometric kernel, where $H$ is the intersection of all index three subgroups of
$\pi_1(F)$.  We have been unable to find reasonable bounds on the size of this
quotient.
\vskip1ex

\section{ A small genus 2 example}

The groups constructed in the previous section are of very large order. If $ F$ is
of genus 2, the corresponding group is of order $2^{38}$. This section presents a
description of a group of order 32, $G_2$, and a homomorphism $\psi_2$ of the
fundamental group of the genus 2 surface to $G_2$, such that $\psi_2$ has a
nongeometric kernel. The next section contains a proof that this example is minimal.

For the remainder of this section $F$ will denote a genus 2 surface.

\vskip3ex

 \noindent{\bf Construction of $G_2$}  For our purposes, the easiest way to describe
$G_2$  is as follows. Define a group structure on the set $(Z_2) ^4 \times  Z_2$
by defining the product by
 $$ (a_1,b_1,a_2, b_2,\epsilon)(a_1',b_1' a_2',b_2',\epsilon') = (a_1
+a_1',b_1+b_1',   b_2+a_2',b_2 +b_2' , \epsilon +\epsilon'+ b_1a_1' +b_2a_2').$$
The operations within the parenthesis are addition and
multiplication in $Z_2$.  The verification that this defines a group
structure can be done by a direct calculation, which is left to the
reader. The group is denoted $32_{42}$  in [6], and $\Gamma_5 a_1$ in the
notation of [3].

An essential calculation for later purposes is that of
commutators in $G_2$.  A direct computation yields
$$[(a_1,b_1,a_2, b_2,\epsilon),(a_1',b_1' a_2',b_2',\epsilon')] =
(0,0,0,0,(b_1a_1' - b_1' a_1) +(b_2a_2'-b_2'a_2)). \eqno{(1)}$$
From this it is apparent
that both the center and commutator subgroup of $G_2$ consists of the set $(0,0,0,0)
\times Z_2$.    The
abelianization of $G_2$ is $(Z_2)^4$, given by the projection $(Z_2)^4 \times Z_2
\rightarrow (Z_2)^4 \times \{0\}$.

\vskip2ex

\noindent{\bf Construction of $\psi_2$ }   Let 
$\{x_1,y_1,x_2,y_2\}$ be a standard generating set of $\pi_1 (F)$   so
that $\pi_1(F)$ has presentation
$< x_1,y_1,x_2, y_2, [x_1,y_1][x_2,y_2 ] =1 >$. This set 
projects to a standard symplectic basis of $H_1 (F;Z_2)$,  $\{
|x_1|,|y_1|,|x_2|,|y_2| \}$

\vskip2ex

		Define $\psi_2 \!: \pi_1(F) \rightarrow G_2$	be setting: 
$$\psi_2(x_1) = (1,0,0,0)\times (0),$$	
$$\psi_2(y_1) = (0,1,0,0)\times (0),$$	
$$\psi_2(x_2) = (0,0,1,0)\times (0), $$	
$$\psi_2(x_2) = (0,0,0,1)\times (0). $$
Using (1) it is easily verified that this gives a well defined
surjective representation. 

The key observation is that $\psi_2$ has the
following property: if $\omega_1$ and $\omega_2$ are elements of $\pi_1(F)$, then

$$[\psi_2(\omega_1), \psi_2(\omega_2)] = (0,0,0,0) \times (|\omega_1| \cap 
|\omega_2|), \eqno{(2)}$$
 where $|\omega_1| \cap 
|\omega_2|$	is the $Z_2$ intersection number of the classes
in $H_1(F;Z_2)$,  represented by $\omega_1$ and $\omega_2$. 	This follows from
(1)  along with the fact that if $(a_1,b_1,a_2, b_2)$ and $(a_1',b_1',a_2', b_2')$
are classes in $H _1(F;Z_2)$, then $(a_1,b_1,a_2, b_2) \cap (a_1',b_1',a_2',
b_2') = (b_1a_1'-b_1'a_1) +(b_2a_2'-b_2'a_2)$. (Note also that the natural map
$\pi_1(F) \rightarrow H_1(F;Z_2 )$ factors through $G_ 2$  via $\psi_2$.)

\vskip2ex

 \noindent{\bf The kernel of $\psi_2$ is nongeometric} Suppose that there is a
simple loop $\gamma$ representing a nontrivial element $\omega$                     
in  Ker$(\psi_2)$.      Our first observation is that
$\gamma$ can be chosen to be separating. If $\gamma$ is nonseparating, pick a simple
loop
$\gamma'$ meeting $\gamma$ transversely in exactly one point. Let $\omega'$ be the
element of $\pi_1(F)$ represented by $\gamma'$.	Clearly, $[\omega,\omega']$ is in
the kernel of $\psi_2$ and it is represented	by a separating simple loop. 

	Since
$\gamma$ is now assumed to be separating, it bounds a
punctured torus on $ F$. This follows from the classification of surfaces. Hence
$\omega = [\omega_1,\omega_2]$, where $\omega_1$ and $\omega_2$ are represented by
simple loops meeting transversely in one point. From this one computes using (2)
that $\psi_2(w) = [\psi_2(\omega_1),\psi_2(\omega_2)] = (0 , 0 , 0 , 0) \times
(|\omega_1| \cap 
|\omega_2|) = (0,0,0,0) \times (1)$, which is nontrivial
in $G_2$. This contradicts the assumption that $w \in$ Ker$(\psi_2).$

\section{ Minimality of G}

The goal of this section is to prove that if $F$ is of genus 2 and the order of $G$
is less than 32, than any homomorphism $\phi \!: \pi_1(F) \rightarrow G$ has geometric
kernel.

 Here is a summary of the argument. We first prove that  any $\phi \!: \pi_1(F)
\rightarrow G$ has geometric kernel  if $G$ is a cyclic extension of an abelian
group, that is, if $G$ contains a normal abelian subgroup with cyclic quotient. The
approach used to prove this was pointed out  by Allan Edmonds. The argument depends
on an analysis of the action of the homeomorphism group of $F$ on the set of
representations of
$\pi_1(F)$ to $G$. We next note that with the exception of two groups of order 24,
$SL_2 (Z_3)$ and
$S_4$, all groups of order less than 32 are cyclic extensions of abelian groups.
This can be proved by a case--by--case analysis using
the Sylow theorems. More easily, group tables such as [6] provide the necessary
information.  The proof is completed using specialized arguments for $SL_2 (Z_3)$
and
$S_4$.

\vskip2ex

\noindent{\bf Cyclic Extension of Abelian Groups}  Fix a group $G$.  The
group of basepoint preserving homeomorphisms of $F$ acts on the set of
representations of $\pi_1(F)$ to $G$, as follows. If $h$ is a homeomorphism of
$F$, send a representation $\phi$ to $\phi \circ h_*$. Notice that $\phi$ has
geometric kernel if and only if $\phi \circ h_*$ has geometric kernel. The following
is a result of Nielsen [4]; a proof can be found in [1].

\beginth{Lemma} If $G$ is a cyclic
group and $\phi \!: \pi(F) \rightarrow G$ is a surjective homomorphism, then there is
a homeomorphism $h$ of $F$ such that $\phi \circ h_*(x_1)$ generates $G$ , and
$\phi \circ h_*(y_1)$, $\phi \circ h_*(x_2)$ and $\phi \circ h_*(y_2)$ are all
trivial. \endth

\beginth{Theorem} If $G$ contains an abelian normal subgroup $N$ such that $G / N$
is cyclic, than any surjective homomorphism $\phi \!: \pi_1(F) \rightarrow G$ has a
geometric kernel. \endth

\noindent {\bf Proof} Denote the quotient map $G \rightarrow G / N$ by
$\rho$.      Applying the lemma, we can assume that $\rho \circ \phi(x_2)$ and
$\rho \circ \phi(y_2)$  are both trivial. Hence $ \phi(x_2)$ and $\rho \circ
\phi(y_2)$ are both in $N$.        The commutator $[x_2,y_2 ]$ is represented by a
simple loop and  is in the kernel of $\phi$, since $\phi([x_2,y_2])$   is in the
commutator subgroup of an abelian group. 
\vskip3ex
\noindent{\bf Exceptional Groups}
\vskip1ex

\noindent{\bf Case 1} We begin by recalling that $SL_2(Z_3)$ is isomorphic to the
semidirect product of the quaternionic 8--group, $Q$, with $Z_3$.    We will use the
standard notation for elements in $Q$ . The generator of $Z_3$  will be denoted
$t $.  The action of $Z_3 $ on $Q$ is given by $tit^{- 1} = j$, $tjt^{- 1} = k$
and $tkt^{-1} = i$.	Note that $t(-1)t^{-1} =-1$  and that $-1 \in Q$ is hence
central in $SL_2 (Z_3)$.

 	Suppose $\phi \!: \pi_1(F) \rightarrow SL_2 (Z_3)$	is a surjective representation
with nongeometric    kernel. Applying the lemma to the composition
$\pi_1(F) \rightarrow SL_2(Z_3) \rightarrow  SL_2(Z_3) / Q = Z_3$ we can assume
that $\phi(x_1) = t q_1$, $\phi(y_1) = q_2$, $\phi(x_2) = q_3$, and $\phi(y_2) =
q_4$, where each $q_i$ is in $Q$.

Since $[x_2 , y_2]$	is represented by a simple loop, $[q_3 , q_4] \ne  1$.
Hence  $[q_3 , q_4] = -1  \in Q$.  It follows that $[tq_1 , q_2] = -1$. Note
that $q_2 \ne \pm 1$,      so $q_2 = \pm i$, $\pm j$, or $\pm k$. From
the commutator relation, $tq_1 q_2 q_1^{-1}t^{-1} = q_2^{-1}$     For any two
quaternions,
$q_1 q_2 q_1^{-1} = q_2^{\pm 1}$.	Hence, $t q_2	t^{-1} = 		q_2^{\pm 1}$.	However,
this is impossible, given that $q_2 \ne \pm 1$  and the action of $t$ on $Q$.

\vskip3ex

\noindent{\bf Case 2} The symmetric group $S_4$ is the semidirect product
of $ Z _2  \times Z_2$ with $S_3$.  As a subgroup, the $Z_2  \times Z_2$ is given by
the set $\{(1),(12)(34),(13)(24),(14)(23)\}$  The $S_3$  is given by the
set $\{(1),(12),(13),(23),(123),(321)\}$.

Let $\phi \!: \pi_1(F) \rightarrow S_4$ be a surjective representation with
nongeometric kernel. The main result of [2] applied to the
composition $\pi_1(F) \rightarrow S_4 \rightarrow S_4 /(Z_2 \times Z_2) \cong S_3$
shows that by applying a homeomorphism we can arrange that $\phi$ takes on the
values $\phi(x_1) = (12)n_1$,   $\phi(y_1) = n_2$,  $\phi(x_2) = (123)n_3$, and 
$\phi(y_2) = n_4$ where each $n_i$   is in $Z_2 \times Z_2$.

     Since both $y_1$ and $y_2$  are represented by simple loops,
neither $n_2$   nor $n_4 $ are trivial.  Also, $[x_1,y_1]$ is represented
by a simple loop, so $[(12)n_1,n_2] \ne 1$. It follows that
$n_2 \ne (12)(34)$. There are two other possibilities for $n_2$.

Suppose that $n_2 = (13)(24)$. There are three possible values
of $n_4$  to be considered. Because $y_1 y_2$ is realized by a simple
loop, $n_4 \ne (13)(24)$.  It is easily seen that $y_1 x_2 y_2^{-1} x_2^{-1}$ 
is realized by a simple loop.  Hence $n_4 \ne (12)(34)$. Finally
$n_4 \ne (14)(23)$, because  $y_2x_1 y_1^{-1} x_1^{-1}$ 
can also be represented by a simple loop.  

	We proceed similarly if $n_2 = 	(14)(23)$.
Clearly  $n_4 \ne (14)(23)$. Because $y_2 x_1y_1^{-1}x_1^{-1}$  is realized by a
simple loop, $n_4 \ne 	 (13)(24)$. Finally, it is again easily seen that 
$y_1 x_2^{-1} y_2^{-1} x_2$    is realized by a simple loop. This implies that $n_4
\ne  (12)(34)$. All possibilities have now been eliminated.

\section{Generalizations}
     The group constructed in Section 3, $G_2$, can be generalized to a group $G_k$
such that for the genus $k$ surface $F_k$ there is a homomorphism $\phi_k \!:
\pi_1(F_k) \rightarrow G_k$  with nongeometric kernel. The arguments are similar
to those of Section 3 and are only outlined here.

     Define $G_k$ by defining a product on the set $(Z_k)^{2k} \times Z_k$ as
follows.
 $$(a_1,b_1, a_2, b_2, \ldots b_k , \epsilon)(a_1',b_1', a_2', b_2',
\ldots b_k' , \epsilon') =(a_1 +a_1',b_1 +b_1', a_2 +a_2', b_2 +b_2', \ldots , b_k
+b_k' , \epsilon +\epsilon' + \sum{b_ia_i'})$$                   
Sums and products within the parenthesis are in $Z_k$.	That this defines a group is
a straightforward calculation. 	

There is a natural representation $\phi_k \!: \pi_1(F_k) \rightarrow G_k$ as before.
In this case the essential observation is

$$   [\phi_k(\omega_1) , \phi_k(\omega_2)] = (0,0,\ldots, 0) \times (|\omega_1|
\cap |\omega_2|),  \eqno{(1)}$$
where $(|\omega_1| \cap |\omega_2|)$ is the $Z_k$ intersection number of
the classes in $H_1(F_k;Z_k)$ represented by $\omega_1$ and $\omega_2$.

If $\phi_k$ had geometric kernel, there would be a separating simple loop
representing an element in the kernel. Using the classification of surfaces, that
element would be of the form
$$[\omega_1 , \omega_1'][\omega_2 , \omega_2']\cdots[\omega_m , \omega_m']$$
with $m < k$ and $(|\omega_i| \cap |\omega_i'|) = 1$   for
all $i$.   A contradiction follows  as in Section 3.

\vskip1ex

\noindent{ \bf Remark} The order of the group $G$ just constructed is $g^{2g+1}$.
This number should be contrasted to the order found in Section 2, $2^{(g-1)2^{ 2g+1}
+2+2g}$. The first is obviously much smaller than the second.  The results of this
paper, along with our difficulties in trying to find smaller examples, leads us to
conjecture that $g^{2g+1}$ represents the least possible order.

\vskip5ex
\centerline{\bf      References}
\vskip2ex

\item{[1]} Edmonds, A. {\sl Surface Symmetry I}, Michigan J. Math 29 (1982)
171--183.

\item{[2]} Edmonds, A. {\sl Surface Symmetry II}, Michigan J. Math 30 (1983)
143--154.

\item{[3]} Hall, M. and Senior, J.K. {\sl The groups of order $2^n (n \le
6)$}     Macmillan, New York (1964).

\item{[4]} Nielson, J. {\sl Die Struktur periodischer Transformationen von 
Flachen}, Dansk Vid. Selsk., Mat.-Fys. Medd. 15 (1937), 1--77.

\item{[5]} Skora, R. {\sl Dissertation}, Department of Mathematics,  University
of Texas, Austin, Texas, 1984.

\item{[6]} Thomas, A.D. and Wood, G.V. {\sl Group Tables}, Shiva  Publishing
Limited, Kent, Great Britian, 1980.

\vskip8ex

\noindent Department of Mathematics

\noindent Indiana University

\noindent Bloomington, IN 47405

\noindent livingst@indiana.edu

\end